\documentclass[twocolumn,aps,floatfix,superscriptaddress,longbibliography]{revtex4-1}
\pdfoutput=1 
\usepackage{amsmath,amssymb,eucal,graphicx,float,epstopdf}
\usepackage{epsfig}
\usepackage[utf8]{inputenc}

\usepackage[colorlinks=true, urlcolor=blue, anchorcolor=blue, citecolor=blue,filecolor=blue,linkcolor=blue,menucolor=blue]{hyperref}

\usepackage{subfigure}

\begin{document}
\title{Dynamic Space Filling}

\author{P.~L.~Krapivsky}
\affiliation{Department of Physics, Boston University, Boston, Massachusetts 02215, USA}
\affiliation{Santa Fe Institute, Santa Fe, New Mexico 87501, USA}

\begin{abstract}
Dynamic space filling (DSF) is a stochastic process defined on any connected graph. Each vertex can host an arbitrary number of particles forming a pile, with every arriving particle landing on the top of the pile. Particles in a pile, except for the particle at the bottom, can hop to neighboring vertices. Eligible particles hop independently and stochastically, with the overall hopping rate set to unity. When the number of vertices in a graph is equal to the total number of particles, the evolution stops when a single particle occupies every vertex. We determine the halting time distribution on complete graphs. Using the mapping of the DSF into a two-species annihilation process, we argue that on $ d$-dimensional tori with $N\gg 1$ vertices, the average halting time scales with the number of vertices as $N^{4/d}$ when $d\leq 4$ and as $N$ when $d>4$.
\end{abstract}

\maketitle

\section{Introduction}

Filling space with identical objects is a fascinating subject \cite{Aste,TS10}.  Two popular realizations are packings and coverings by balls. Balls cannot overlaps in a packing, while they do overlap in coverings as each point should be covered (i.e., belong to at least one ball). The densest sphere packings and the least dense sphere coverings of $\mathbb{R}^d$ are particularly popular subjects \cite{Kershner39,Rogers,Conway,Hales05,Vallentin05,Vallentin06, Hales11,Cohn16,Vallentin16,Cohn17}. The densest sphere packings are known \cite{Hales11,Cohn16,Vallentin16,Cohn17} for $d=1,2,3,8, 24$. The least dense sphere coverings are known \cite{Kershner39,Vallentin05,Vallentin06} for $d=1,2$ and conjecturally for $d=24$. 

Random sequential adsorption (RSA) is the dynamic counterpart of packing where deposition events leading to the overlap with already present balls are discarded \cite{Evans93, Talbot00, Tor, KRB}. Random sequential covering (RSC) is a dynamic counterpart of covering where deposition events leading to an increase of coverage are accepted \cite{Krapivsky23, Pascal}. In RSA and RSC processes, the evolution stops when the system reaches a jammed state. Jammed states are random, albeit the final jamming density of the infinite space $\mathbb{R}^d$ is deterministic. The RSA process stops at the filling fraction $\theta_\text{RSA}(d)<1$ depending on the spatial dimension; the RSC process stops at the covering number (the average number of balls covering a point) $\theta_\text{RSC}(d)>1$. The filling and covering numbers are analytically known only in one dimension. 

We investigate the dynamic space filling (DSF) leading to perfect space filling. Our `space' is a graph with particles occupying vertices; the perfect filling is the arrangement with a single particle at any vertex. The DSF proceeds independently on the connected components of a graph. Therefore, we can limit ourselves to connected graphs. To avoid unnatural behaviors, we consider simple graphs  \cite{Diestel}, i.e., undirected graphs without self-loops and multiple edges. We also limit ourselves to regular graphs \cite{Diestel} in which each vertex has the same number of neighbors. Simple connected $r-$regular graphs with the lowest or highest number $r$ of neighbors are easy to classify: The only 2-regular connected graph with $N$ vertices is the ring $R_N$, and the only $(N-1)-$regular graph with $N$ vertices is the complete graph $K_N$. For $2<r<N-1$, the total number of $ r$-regular graphs is unknown (rapidly growing with $N$ for any fixed $r$). 

The DSF process on a graph with $N$ vertices concerns the dynamics of $N$ random walkers. Random walkers form a pile at each vertex according to the order of their arrival at the vertex: The first arrival is at the bottom, etc. The random walker at the bottom is stuck to the vertex, while other random walkers hop independently and stochastically to neighboring vertices. We set the overall hopping rate to unity, so the hopping rate to each neighbor is $1/r$ in the case of $r-$regular graphs. The DSF process comes to a halt when each vertex is occupied by exactly one random walker. 

At first sight, the DSF process looks significantly simpler than the lattice versions of RSA and RSC, such as packing or covering of $\mathbb{Z}^d$ with dimers. Indeed, for the RSA and RSC, the final jammed states are understood only in one dimension; for the DSF process, the jammed state is a perfectly filled graph, i.e., it is universal and trivial. However, the approach to the perfect filling is highly non-trivial already in one dimension.

Although our terminology has a mathematical and computer science flavor (regular graphs, halting time distributions, etc.), we emphasize that the DSF process admits a mapping into a diffusion-controlled annihilation process. Indeed, the DSF process on a connected graph is isomorphic to the two-species diffusion-controlled annihilation process on the same graph. We denote the particles of the two species by $A$ and $B$ and map the DSF representation into the $A-B$ representation as follows:
\begin{itemize}
\item An empty vertex in the DSF hosts a $B$ particle.
\item A vertex occupied by a single particle in the DSF is empty in the realm of the annihilation process. 
\item A vertex occupied by $k+1$ particles in the DSF hosts $k$ particles of type $A$. 
\end{itemize}
Therefore, $A$ and $B$ particles never share the same vertex: Each vertex is either empty or occupied, and in the latter case, it either contains a single $B$ particle or an arbitrary number of $A$ particles. The rules of the DSF process imply that non-interacting $A$ particles undergo independent, identical random walks. We call $A$ particles active since they are mobile; passive $B$ particles are immobile.

When an $A$ particle hops into an empty vertex, i.e., the vertex occupied by $B$ particle, the vertex gets perfectly filled, i.e., $A$ and $B$ particles immediately annihilate:
\begin{equation}
\label{AB}
A+B\to \emptyset
\end{equation}
In other words, we have a two-species diffusion-controlled annihilation process. We postulate that the number of particles in the DSF process is equal to the number of vertices so that the perfect filling can be achieved. In terms of the annihilation process, this means that the initial numbers of $A$ and $B$ particles are equal; then, the numbers of $A$ and $B$ particles will remain equal forever and the perfect filling corresponds to the vacuum state of the annihilation process. 

The annihilation process \eqref{AB} with one species immobile mimics several chemical processes such as dissolution of solids, corrosion, and etching, see \cite{Meakin-IDLA,dissol-IDLA, KM-IDLA, Havlin-IDLA}. The initial conditions appropriate for these chemical processes differ from the ones that ensure the perfect filling. It is customary to start with an infinite lattice with each site occupied by an immobile $B$ particle, and at time $t=0$ turn on a localized source injecting $A$ particles into a small localized region. The $B$ particles disappear inside a growing droplet, and the relevant questions concern the size and the shape of the droplet \cite{Havlin-IDLA}, the survival probability of an $A$ particle long after it was injected \cite{PK-IDLA}, etc. Starting instead with all $A$ particles in a small localized region, one can deduce heuristic predictions for the DSF process on $d$-dimensional tori subject to the localized initial conditions. 

The rest of this paper is organized as follows. In Sec.~\ref{sec:complete}, we analyze the DSF process on complete graphs. We determine the halting time distribution for any $N$. The halting time increases linearly with $N$ and remains a non-self-averaging random variable. In the $N\to\infty$ limit, the halting time distribution approaches the scaling form: $P(T_N)\to N^{-1}\mathcal{P}(\tau)$ with $\tau=T_N/N$. The scaled distribution reads
\begin{equation}
\label{P-tau}
\mathcal{P}(\tau) = \sum_{k=-\infty}^\infty  (-1)^{k+1}\, k^2\, e^{-k^2\tau} 
\end{equation}
Alternatively, $\mathcal{P}(\tau)$ can be expressed via a theta function.

We also briefly discuss the DSF on tori $T_d(L)=(R_L)^d$, equivalently hypercubes with periodic boundary conditions (Sec.~\ref{sec:d-dim}). Such tori are regular graphs with $N=L^d$ and $r=2d$. The DSF process on the infinite lattices $\mathbb{Z}^d$ with a density of random walkers equal to unity relaxes to the perfectly filled state. We use the mapping into a two-species annihilation process with equal concentrations of both species and one species immobile. The vacuum state of this annihilation process corresponds to the perfect filling of the DSF process. Relying on the asymptotic decay laws for the densities in the two-species annihilation process on the infinite hyper-cubic lattices $\mathbb{Z}^d$ we estimate a typical halting time on the tori with $N$ vertices:
\begin{equation}
\label{T-d}
T_N\sim
\begin{cases}
N^{4/d}        & d < 4  \\
N                & d \geq 4
\end{cases}
\end{equation}
These estimates are conjectural since the arguments leading to \eqref{T-d} are heuristic. The nature of the random variable $T_N$ on large tori is unsettled when $d<4$; for $d\geq 4$, the halting time $T_N$ is a non-self-averaging random variable as we argue in Sec.~\ref{sec:d-dim}.  

Apart from the halting time $T_N$, the duration $t_N$ of the last step, namely, the evolution from a state with $N-2$ perfectly filled vertices to the final perfect filling, is another interesting random quantity. In Sec.~\ref{sec:disc}, we outline how to compute the joint distribution $P(T_N, t_N)$ on the complete graphs. 

In computing the halting time distribution $P(T_N)$ and the joint distribution $P(T_N, t_N)$ on the complete graphs, we rely on the fact that the evolution of the total number $m$ of empty vertices, $m\to m-1$,  is determined by $m$ alone. This crucial property occurs only on complete graphs; generally, the neighbors of each empty vertex and their occupancies matter. One can compute the distributions $P(T_N)$ and $P(T_N, t_N)$, circumventing the detailed knowledge of the evolution of the random variable $m$. However, our methods give the exact Laplace transform of the probability distribution $P_m(t)$. 

In Appendix~\ref{ap:der}, we compute the average and the variance of $m$ on large complete graphs and show that $P_m(t)$ becomes Gaussian when $N\gg 1$. In Appendix~\ref{ap:high}, we show how to determine higher cumulants of $m$. The computations quickly become cumbersome, so we have computed only the third and the fourth cumulants. Intriguingly, the same values of ratios of the cumulants to the average have appeared in the context of the statistics of current in quantum conductors and symmetric exclusion process \cite{Levitov,Blanter,Nagaev,Derrida04,Bodineau04,Bodineau13}. In Appendix~\ref{ap:loc}, we analyze the annihilation process \eqref{AB} on the inifinite hypercubic lattices with localized initial condition, and use the results to guess the scaling of the halting time of the DSF process on tori with locatlized initial conditions.

\section{DSF on complete graph}
\label{sec:complete}

Consider the DSF process on complete graphs. It is convenient to take the complete graph $K_{N+1}$ so that each active particle hops with rate $1/N$ to every of $N$ neighboring vertices. (Recall that we set the overall hoping rate to unity.) Denote by $[m]$ the state with $m$ empty vertices in the DSF representation. The transition $[m]\to [m-1]$ occurs with rate $r_m=m^2/N$ as there are $m$ active $A$ particles and $m$ passive $B$ particles. If $t_m$ is the transition time from $m$ to $m-1$,  the halting time is $T_N=t_1+t_2+\cdots+t_{m_0}$, where $m_0$ is the initial number of empty vertices. The transition times are exponentially distributed:
\begin{equation}
\label{tm}
\Pi(t_m)=r_m\,e^{-r_m t_m}
\end{equation}
From \eqref{tm}, the average transition time is $\langle t_m\rangle=r_m^{-1}$ and the variance is $\langle\!\langle t_m^2\rangle\!\rangle=\langle t_m^2\rangle-\langle t_m\rangle^2=r_m^{-2}$. Therefore the average halting time is
\begin{equation}
\label{TN-av}
\langle T_N\rangle=\sum_{m=1}^{m_0} \langle t_m\rangle=N\sum_{m=1}^{m_0} m^{-2}
\end{equation}
In the $m_0\to\infty$ limit, the sum converges to the special value of the zeta function \cite{NT-1} discovered by Euler: $\sum_{m\geq 1} m^{-2}=\zeta(2)=\frac{\pi^2}{6}$. Therefore in the leading order
\begin{equation}
\label{TN}
\langle T_N\rangle=\frac{\pi^2}{6}\,N
\end{equation}
The second moment $\langle T_N^2\rangle=\langle T_N\rangle^2 +\sum_{1\leq m\leq m_0} \langle\!\langle t_m^2\rangle\!\rangle$ becomes
\begin{equation}
\label{2TN-av}
\langle T_N^2\rangle=\langle T_N\rangle^2+N^2\sum_{m=1}^{m_0} m^{-4}
\end{equation}
Recalling $\sum_{m\geq 1} m^{-4}=\zeta(4)=\frac{\pi^4}{90}$ also discovered by Euler we obtain
\begin{subequations}
\begin{equation}
\label{mu-2}
\mu_2=\lim_{N\to\infty}\frac{\langle T_N^2\rangle}{\langle T_N\rangle^2} =  \frac{7}{5}
\end{equation}
The initial value $m_0\sim N$ and hence replacing the sums in \eqref{TN-av} and \eqref{2TN-av} by infinite sums is exact up to $N^{-1}$ in the first case and up to $N^{-3}$ in the second. The precise value  value of $m_0$ does not affect the leading behavior. (The random initial distribution gives $m_0\approx N/e$, see \eqref{IC}. But the same asymptotic behaviors emerge for the extreme initial condition $m_0=N$ describing the evolution starting from all particles initially at a single vertex.)

The normalized moments 
\begin{equation*}
\mu_p=\lim_{N\to\infty}\frac{\langle T_N^p\rangle}{\langle T_N\rangle^p}
\end{equation*}
appear to be rational for all integer positive $p$. Laborious but straightforward calculations give
\begin{equation}
\label{mu-345}
\mu_3 =  \frac{93}{35}\,, \quad  \mu_4 =  \frac{1143}{175}\,, \quad \mu_5 =  \frac{219}{11}
\end{equation}
\end{subequations}
The rationality of all moments $\mu_p$ is proven below. 

The moments are non-trivial implying that the halting time is an asymptotically non-self-averaging random variable. Thus for its complete characterization, we must determine the halting time distribution. The halting time distribution $P(T_N)$ can be extracted from the probability distribution $P_m(t)$ of the state of the system. This probability distribution obeys 
\begin{equation}
\label{Pm}
\frac{d P_m}{dt}=r_{m+1}P_{m+1}-r_mP_m
\end{equation}
and the initial condition  $P_m(0)=\delta_{m,m_0}$. The halting time distribution is then found from $P(T_N)=r_1P_1(T_N)$. To treat Eqs.~\eqref{Pm} we use the Laplace transform: 
\begin{equation}
\label{Qms:def}
Q_m(s)=\int_0^\infty dt\,e^{-st}\,P_m(t)
\end{equation}
The Laplace transform of the halting time distribution is
\begin{eqnarray}
\label{Lap-halt}
Q(s)=\int_0^\infty dt\,e^{-sT_N}\,P(T_N)=r_1Q_1(s). 
\end{eqnarray}
Performing the Laplace transform of Eqs.~\eqref{Pm} yields  
\begin{subequations}
\label{Qms-rec}
\begin{equation}
\label{Qms-eq}
(s+r_m)Q_m(s)=r_{m+1}Q_{m+1}(s)
\end{equation}
for $1\leq m\leq m_0-1$ and
\begin{equation}
\label{Qms-0}
(s+r_{m_0})Q_{m_0}(s)=1
\end{equation}
\end{subequations}
Starting with \eqref{Qms-0} and iterating \eqref{Qms-eq} we determine $Q_m(s)$ for all $m$. In particular, the Laplace transform \eqref{Lap-halt} of the halting time distribution reads 
\begin{equation}
\label{Qs}
Q(s)=\prod_{m=1}^{m_0} \frac{r_m}{r_m+s}
=\prod_{m=1}^{m_0} \left[1+\frac{sN}{m^2}\right]^{-1}
\end{equation}

When $N\to\infty$, the halting time distribution approaches the scaling form
\begin{equation}
\label{PTN}
P(T_N)\to N^{-1}\,\mathcal{P}(\tau) \quad {\rm with}\quad  \tau=\frac{T_N}{N}
\end{equation}
implying that the Laplace transform admits the scaling
form
\begin{equation}
\label{Qss}
Q(s)\to \mathcal{Q}(\sigma), \qquad \sigma=s N
\end{equation}
in the $N\to\infty$ limit.  This indeed agrees with Eq.~\eqref{Qs} that also gives the scaling form of the Laplace transform
\begin{equation}
\label{Qsigma}
\mathcal{Q}(\sigma)
=\prod_{m\geq 1} \left[1+\frac{\sigma}{m^2}\right]^{-1}=\frac{\pi\sqrt{\sigma}}{\sinh\big(\pi\sqrt{\sigma}\big)}
\end{equation}

The expansion of $\mathcal{Q}(\sigma)$ near $\sigma=0$ confirms previous results \eqref{mu-2} and \eqref{mu-345} found after laborious calculations. Our empirical observation about the rationality of $\mu_p$ also becomes obvious. Expanding $\mathcal{Q}(\sigma)$ to higher orders, one can determine any desirable $\mu_p$. For instance,
\begin{equation*}
\mu_6 =  \frac{12730293}{175175}\,, \quad  \mu_7 =  \frac{221157}{715}\,, \quad \mu_8 =  \frac{457141779}{303875}
\end{equation*}

\begin{figure}
\centering
\includegraphics[width=7.89cm]{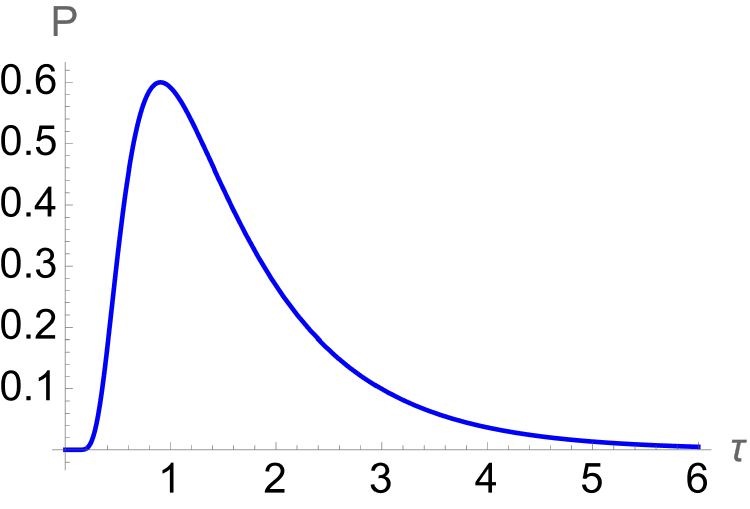}
\caption{The scaled halting time distribution $\mathcal{P}(\tau)$.  }
\label{Fig:Ptau}
\end{figure}

Inverting Laplace transform \eqref{Qsigma} yields the announced distribution \eqref{P-tau} of the scaled halting time. The asymptotic behaviors of the scaled halting time distribution are 
\begin{equation}
\label{P-tau-asymp}
\mathcal{P}(\tau) \simeq
\begin{cases}
\frac{1}{2}\big(\frac{\pi}{\tau}\big)^\frac{5}{2} e^{-\frac{\pi^2}{4\tau}}  & \tau\to 0\\
2 e^{-\tau}  & \tau \to \infty
\end{cases}
\end{equation}
(see also Fig.~\ref{Fig:Ptau}). The large time asymptotic follows from the exact solution \eqref{P-tau}: When $\tau\to\infty$,  the sum on the right-hand side of Eq.~\eqref{P-tau} is dominated by the terms with $k=\pm 1$. The small time asymptotic is extracted from the $\sigma\to\infty$ behavior of the Laplace transform \eqref{Qsigma}.

During the initial stage far before the halting time, viz. when $t\ll N$, the number of $B$ particles is an asymptotically self-averaging random quantity. The average 
\begin{subequations}
\label{m-av-var}
\begin{align}
\label{m-av}
\langle m\rangle = N n, \qquad n=\frac{1}{e+t}
\end{align}
is proportional to $N$. The variance $\langle\!\langle m^2\rangle\!\rangle \equiv \langle m^2\rangle-\langle m\rangle^2$ is also proportional to $N$ for large $N$:
\begin{align}
\label{m-var}
\langle\!\langle m^2\rangle\!\rangle = N v, ~\quad v=\frac{1}{3}\left[\frac{1}{e+t} + \frac{2e^3-3e^2}{(e+t)^4}\right]
\end{align}
\end{subequations}
These formulae for the density $n(t)$ and the rescaled variance $v(t)$ are valid in the case when initial positions of the particles in the DSF process are uncorrelated. Indeed, the derivation of Eqs.~\eqref{m-av-var}, see Appendix~\ref{ap:der}, takes into account that for the annihilation process $n(0)=e^{-1}$ and $v(0)=e^{-1}(1-e^{-1})$ if the initial positions of the particles in the DSF process are uncorrelated and $N\gg 1$. 

The distribution $P_m(t)$ acquires the scaling form
\begin{subequations}
\label{Pmt-Phi}
\begin{equation}
\label{Pmt-scal}
P_m(t)=\frac{1}{\sqrt{Nv}}\,\Phi(\xi)
\end{equation}
in the scaling region
\begin{equation}
\label{scal:def}
N\to\infty, \quad t\to\infty, \quad \xi = \frac{m-Nn}{\sqrt{Nv}}=\text{finite}
\end{equation}
The scaled distribution is Gaussian with zero mean and unit variance
\begin{equation}
\label{Gauss}
\Phi(\xi) = \frac{1}{\sqrt{2\pi}}\,e^{-\xi^2/2}
\end{equation}
\end{subequations}
The proof of \eqref{m-av-var}--\eqref{Pmt-Phi} is relegated to the Appendix~\ref{ap:der}. 

Denote by $F(t_N)$ the distribution of the duration $t_N$ of the final step. This distribution can be found from the same formula $F(t_N)=r_1P_1(t_N)=N^{-1}P_1(t_N)$ as the distribution of the halting time, the difference is that  as the initial condition we should use $P_m(0)=\delta_{m,1}$ asserting that we start with a single empty vertex. Solving
\begin{equation}
\label{P1}
\frac{d P_1}{dt_N}=-\frac{1}{N}\,P_1, \qquad P_1(0)=1
\end{equation}
yields $F(t_N)=N^{-1}e^{-\tau}$ with $\tau=t_N/N$. Thus, the duration of the last step is a non-self-averaging random variable. The normalized moments are
\begin{equation}
\label{ratios}
\frac{\langle t^p\rangle}{\langle t\rangle^p}=p!
\end{equation}
The exponential distribution of the duration of the final time is exact for arbitrary $N$, and hence Eqs.~\eqref{ratios} are exact and independent on $N$. 

For the localized initial condition (all particles in the DSF process are initially at a single vertex), the initial conditions for the annihilation process are $n(0)=1$ and $v(0)=0$ when $N\gg 1$, and instead of \eqref{m-av-var} we obtain
\begin{equation}
\label{n-av-var}
n=\frac{1}{1+t}\,, \qquad v=\frac{1}{3}\left[\frac{1}{1+t} - \frac{1}{(1+t)^4}\right]
\end{equation}
as we show in Appendix~\ref{ap:der}.

\section{DSF in finite dimensions}
\label{sec:d-dim}

We begin with the DSF process on the infinite hyper-cubic lattices $\mathbb{Z}^d$. The evolution continues forever. We are mostly interested in the long time behavior. If the initial positions are uncorrelated, the probability $p_\ell$ to find $\ell$ particles at a site is 
\begin{equation}
\label{IC}
p_\ell = \frac{e^{-1}}{\ell!}
\end{equation}

We map the DSF process onto a two-species diffusion-controlled annihilation process \eqref{AB}. Our convention that the overall hopping rate of $A$ particles is equal to unity implies that the diffusion coefficient of active particles is $1/(2d)$ on the hypercubic lattice $\mathbb{Z}^d$. The initial condition \eqref{IC} for the DSF process leads to $b(0)=p_0$, i.e., the initial densities of $A$ and $B$ particles are
\begin{equation}
\label{ab:IC}
a(0)=b(0)=e^{-1}
\end{equation}
The densities of $A$ and $B$ particles remain equal throughout the evolution, $a(t)=b(t)=n(t)$, and decay as
\begin{equation}
\label{nt-d}
n(t)\sim
\begin{cases}
t^{-d/4}        & d < 4  \\
t^{-1}           & d \geq 4
\end{cases}
\end{equation}
when $t\gg 1$. A naive `mean-field' treatment suggests that the density obeys $\dot n = -n^2$ leading to the $t^{-1}$ decay \cite{KRB}. The decay laws \eqref{nt-d} thus assert that the mean-field treatment provides qualitatively correct description above the critical dimension $d\geq d_c=4$. 

The chief reason for the slowing down, $t^{-d/4}$ versus the mean-field $t^{-1}$ decay, is spontaneous `phase' separation in $d<4$ dimensions: The system organizes \cite{Zeldovich78,Wilczek83,Redner84} into a mosaic of $A$ and $B$ domains. This phenomenon is particularly striking in one dimension with alternating $A$ and $B$ domains and annihilation process \eqref{AB} at the domain boundaries \cite{LR91a,LR91b,LR92}. The one-dimensional mosaic is characterized by the three basic length scales, the diffusive $t^{1/2}$ length of domains of similar particles, the interparticle $t^{1/4}$ spacing, the $t^{3/8}$ depletion zone between adjacent domains \cite{LR91a,LR91b,LR92}, and the mosaic looks like
\begin{equation}
\label{mozaic}
\cdots  \underbrace{A\quad \overbrace{A \quad A}^{t^{1/4}} \cdots A \quad A}_{t^{1/2}} \underbrace{\qquad}_{t^{3/8}}\underbrace{B\quad B \cdots B \quad B}_{t^{1/2}}  \cdots
\end{equation}
The $\ell_{AA}(t)\sim t^{1/4}$ length scale represents the typical separation between adjacent mobile $A$ particles. The average interparticle separation between immobile particles is also $\langle \ell_{BB}(t)\rangle \sim t^{1/4}$, but the density of closest neighbors separated by small distance $x$ remains finite and scales as $x^{-3/2}$ for small $x$,  see \cite{LR92}. Hence, the reduced moments $B_b(t)\equiv \langle \ell^b_{BB}(t)\rangle^{1/b}$ remain finite for $b < \frac{1}{2}$, and grow algebraically, $B_b\sim t^{(2b-1)/4b}$ for $b > \frac{1}{2}$. 

The decay laws \eqref{nt-d} have been established via heuristic \cite{Zeldovich78,Wilczek83,Redner84} and rigorous \cite{Bramson88,Bramson91} analyses for the symmetric version of the two-species annihilation process \eqref{AB} when the diffusion coefficients are equal: $D_A=D_B$. The chief prediction \eqref{nt-d} is expected to hold when both diffusion coefficients are positive: $D_A\geq D_B>0$. Our $B$ particles are immobile. The proof \cite{Bramson88,Bramson91} of the decay laws \eqref{nt-d} assumes that both species are mobile, but it seemingly can be extended to the case when one species is immobile. 

The spatial organization in the annihilation process with one immobile species was studied mostly in one dimension.  The emerging spatial mosaic \cite{LR92} is significantly different from the case when both species have equal diffusivities. Nevertheless, the validity of the decay laws \eqref{nt-d} in the extreme case when one species is immobile is supported by simulations \cite{LR92,Blumen91,Blumen96} in the case of equal concentarations. (Qualitatively different behaviors emerge in the case of unequal concentrations \cite{Blumen91,Blumen96}.) 

Direct numerical simulations of the two-species annihilation is challenging. The mean-field $t^{-1}$ decay holds for a long time in $d=2,3$ before the emergence of the ultimate $t^{-d/4}$ asymptotic. Confirming the $t^{-1/2}$ decay in two dimensions has proven difficult, and the $t^{-3/4}$ decay in three dimensions is essentially impossible to observe. Only the one-dimensional $t^{-1/4}$ decay quickly emerges. A numerical integration of {\em deterministic} partial differential equations with {\em random} intial densities and reaction term $a({\bf r},t)b({\bf r},t)$ is more amenable in two and three dimensions, and the results support the decay laws \eqref{nt-d} for the extreme case of immobile $B$ particles \cite{Blumen91,Blumen96}.  

Combining the decay laws \eqref{nt-d} and criterion 
\begin{equation}
\label{crit}
N n(T)\sim 1
\end{equation}
we arrive at the scaling laws \eqref{T-d} for the halting time $T_N$. The criterion \eqref{crit} is difficult to justify as fluctuations might be strong when the expected number of particles $N n(t)$ is small. To address this issue analytically, it is customary to represent a random variable as a sum of an average that is deterministic and scales linearly with $N$ and a stochastic contribution proportional to $\sqrt{N}$. This so-called Van Kampen expansion \cite{VanKampen} provides an adequate description of many reaction processes, see, e.g., \cite{Spouge85b,Sid86,Ernst87a,Hilhorst04,McKane,BK-van}. In the present situation 
\begin{equation}
\label{van}
N_A=Nn(t) +\sqrt{N}\,\eta
\end{equation}
The stochastic variable $\eta$ has zero mean, $\langle \eta\rangle=0$, and unknown variance $\langle \eta^2\rangle=v(t)$. If the stochastic contribution is sub-dominant, $\sqrt{N v(t)}\lesssim Nn(t)$ during the evolution process, the naive criterion $N n(T)\sim 1$ can be used to estimate the halting time. Otherwise, the criterion $\sqrt{N v(T)}\sim Nn(T)$ is more appropriate. This latter criterion is difficult to justify, and if it works \cite{BK-van,Wendy21,Wendy23,Sergey23}, one should know $v(t)$ to estimate the halting time. The variance $v(t)$ is unknown for the annihilation process \eqref{AB} in $d$ dimensions. (In Appendix~\ref{ap:der}, we compute $v(t)$ on the complete graphs and show that on complete graphs one can use the naive criterion $N n(T)\sim 1$.)

Denote by $(1,1)$ the state before halting (a single active $A$ particle and a single passive $B$ particle).  Eventually $(1,1)\to (0,0)$ and the evolution stops.  Let $\mathcal{T}_N$ be a typical time between the formation of the $(1,1)$ state and its disapperance. This time, equivalently the time for a random walker to hit a stationary target, scales as
\begin{equation}
\label{tau-d}
\mathcal{T}_N \sim
\begin{cases}
N^2        & d = 1  \\
N\ln(N)   & d = 2 \\
N            & d \geq 3
\end{cases}
\end{equation}
Comparing \eqref{T-d} and \eqref{tau-d} we see that $T_N$ greatly exceeds $\mathcal{T}_N$ in dimensions $d=1,2,3$; in higher dimensions, $d\geq 4$, the times $T$ and $\mathcal{T}$ are comparable. (We digress and note that as usuual, the behavior in two dimensions is most subtle \cite{Aldous89,Hilhorst91,Peres04, Peres07,Peres-BM}. Still, a lot is known, e.g., the number of steps it takes for a random walker to completely {\em cover} the torus is a self-averaging random variable growing as $CN[\ln N]^2$ with known amplitude \cite{Peres04}.)

The duration $\mathcal{T}$ of the final step before halting is a non-self-averaging random variable. Since $T_N\sim \mathcal{T}_N$ when $d\geq 4$, the halting time remains a non-self-averaging random variable in the $N\to\infty$ limit when $d\geq 4$. The same probably holds also when $d<4$. To gauge this feature quantitatively, one can compare the average $\langle T_N\rangle$ and the variance $\langle\!\langle T_N^2\rangle\!\rangle = \langle T_N^2\rangle- \langle T_N\rangle^2$. The variance is expected to scale algebraically 
\begin{equation}
\label{Var-T-d}
\langle\!\langle T_N^2\rangle\!\rangle \sim
\begin{cases}
N^{\beta_d}    & d < 4  \\
N^2                & d \geq 4
\end{cases}
\end{equation}
If $\beta_d\geq \frac{4}{d}$ when $d<4$, the halting time is a non-self-averaging random quantity in low spatial dimensions. 

The initial conditions play a significant role in the DSF process in finite dimensions. To illustrate this statement, consider for the DSF process in one dimension starting from the antiferromagnetic, in the $A-B$ particles language, initial condition: $\ldots ABABAB\ldots$. The spatial distribution remains antiferromagnetic throughout the evolution
\begin{equation}
\label{mozaic-alt}
\cdots\quad  A\quad B\quad A\quad B\quad A\quad B\quad \cdots
\end{equation}
Adjacent particles are not necessarily on the neighboring sites. The typical separation is $\langle \ell_{AB}(t)\rangle \sim t^{1/2}$. Thus, $a(t)=b(t)\sim t^{-1/2}$ and the typical halting time scales as $T_N\sim N^2$ with the length $N$ of the ring. The duration of the final step before halting is comparable, $\mathcal{T}_N\sim N^2$. Summarizing, the halting time for the DSF process on the ring starting from the antiferromagnetic initial condition is a non-self-averaging random variable with 
\begin{equation}
\label{av-var-AF}
\langle T_N\rangle\sim N^2, \qquad \langle\!\langle T_N^2\rangle\!\rangle \sim N^4
\end{equation}

For the DSF process on complete graphs, even in the extreme case of the localized initial distribution, the asymptotic behavior is the same as for the uniform initial distribution. To probe the influence of the localized initial distribution for the DSF process in finite dimensions, it is natural to begin with infinite initially empty lattices, except for $N\gg 1$ particles at the origin. Such processes on $\mathbb{Z}^d$ admit a comprehensive analytical description detailed in Appendix \ref{ap:loc}. Combining such analytical results for the infinite systems with heuristic arguments, we estimated the lower bounds of the halting times on tori with $N$ vertices, $T_N > N^\frac{2}{d}$. The lower bound gives an exact asymptotic only in one dimension. We finally note that the inequality $T_N > \mathcal{T}_N$ valid independently on the initial conditions provides a stronger lower bound, cf. \eqref{tau-d}.

\section{Discussion}
\label{sec:disc}

The DSF process on complete graphs (Sec.~\ref{sec:complete}) is solvable. We computed the distribution of the halting time and the distribution of the duration of the last step for arbitrary $N$. The former distribution becomes particularly neat in the $N\to\infty$ limit when it approaches the scaled form \eqref{P-tau}. The same scaled distribution \eqref{P-tau} describes the width distribution of one-dimensional interfaces modeled by periodic random walks \cite{Zoltan94}. The reason for this unexpected coincidence, like a hidden isomorphism between two seemingly very different problems, is unknown. 

In computing the halting time distribution $P(T_N)$, we relied on the fact that the evolution of the total number $m$ of empty vertices, $m\to m-1$,  is determined by $m$ alone. This crucial property occurs only on complete graphs; generally, the neighbors of each empty vertex and their occupancies matter. We computed $P(T_N)$, circumventing the detailed knowledge of the evolution of $m$. The random variable $m$ is interesting in its own right. In Appendices \ref{ap:der} and \ref{ap:high}, we probe the statistics of $m$. The computed ratios of the cumulants up to the fourth order to the average are the same as for the current in quantum conductors and the symmetric exclusion process \cite{Levitov,Blanter,Nagaev,Derrida04,Bodineau04,Bodineau13}. We conjecture that the same equality holds for all cumulants.

One can probe more convoluted temporal characteristics of the DSF process on complete graphs. For instance, the joint distribution $P(T_N, t_N)$ of the duration $T_N$ of the process (the halting time) and the duration $t_N$ of the last step, $t_N<T_N$, approaches the scaling form 
\begin{subequations}
\label{PT2}
\begin{equation}
\label{PTt}
P(T_N, t_N)=N^{-2}\mathcal{R}(\tau-\tau')\,e^{-\tau'}
\end{equation}
with $(\tau,\tau')=N^{-1}(T_N,t_N)$. The Laplace transform of the scaled distribution $\mathcal{R}(\tau)$ reads 
\begin{equation}
\label{Rsigma}
\int_0^\infty d\tau\,e^{-\sigma \tau}\mathcal{R}(\tau)
= \prod_{m\geq 2} \frac{1}{1+\frac{\sigma}{m^2}} = \frac{\pi\sqrt{\sigma}\,(1+\sigma)}{\sinh\big(\pi\sqrt{\sigma}\big)}
\end{equation}
\end{subequations}
When $t_N<T_N\ll N$, the joint distribution \eqref{PT2} becomes
\begin{equation*}
P(T_N, t_N)\simeq \left(\frac{\pi N}{T_N-t_N}\right)^\frac{9}{2}\,\frac{\exp\!\left[-\frac{\pi^2 N}{4(T_N-t_N)}-\frac{t_N}{N}\right]}{8N^2}
\end{equation*}

Needless to say, the variables $T_N$ and $t_N$ are correlated. For instance, from \eqref{PTt} one finds
\begin{subequations}
\begin{equation}
\frac{\langle T_N t_N\rangle}{N^2} = \int_0^\infty d\tau\,(\tau+2)\mathcal{R}(\tau)
\end{equation}
and from \eqref{Rsigma} one computes
\begin{equation}
 \int_0^\infty d\tau\,\mathcal{R}(\tau) = 1, \quad    \int_0^\infty d\tau\,\tau\mathcal{R}(\tau)=\frac{\pi^2}{6}-1
\end{equation}
\end{subequations}
implying that 
\begin{equation}
\frac{\langle T_N t_N\rangle}{N^2} = \frac{\pi^2}{6}+1 > \frac{\pi^2}{6} = \frac{\langle T_N\rangle \langle t_N\rangle}{N^2}
\end{equation}
A similar calculation yields
\begin{equation}
\frac{\langle T_N t_N^p\rangle}{N^{p+1}} = p!\,\frac{\pi^2}{6} + p\cdot p! > p!\,\frac{\pi^2}{6} = \frac{\langle T_N\rangle \langle t_N^p\rangle}{N^{p+1}}
\end{equation}

The DSF process in finite dimensions (Sec.~\ref{sec:d-dim}) admits a mapping into two-species annihilation process with exactly equal numbers of particles of both species and particles of one species performing nearest-neighbor random walk and immobile particles of another species. This mapping leads to the decay laws \eqref{nt-d} for the densities of the species which we combined with naive criterion \eqref{crit} to estimate the halting time \eqref{T-d}. More work is required to justify or disprove \eqref{T-d} in $d<4$. 

Conjecturally, the halting time remains a non-self-averaging random variable even in the $N\to\infty$ limit. One can justify this assertion when $d\geq 4$ by noting that (i) the duration of the final step is a non-self-averaging random variable, and (ii) the duration of the final step is comparable with halting time when $d\geq 4$. If  $d\leq 3$, the duration of the final step is asymptotically negligible compared to the halting time, so we cannot use the same argument as in $d\geq 4$.

An extension to regular random graphs \cite{Wormald} is an interesting avenue for future work. A regular random graph with degree $r\geq 3$ picked uniformly from all possible $r-$regular graphs is effectively infinite-dimensional as manifested by a logarithmic growth of the diameter with $N$. We anticipate that the halting time remains a non-self-averaging random variable in the $N \to \infty$ limit. The distribution of the halting time should approach a scaling form \eqref{PTN}, and only the scaled distribution probably depends on $r$. 

It would be interesting to analyze {\em planar} regular random graphs (see \cite{Francesco95,Francesco02} and references therein). The remarkable feature of planar regular random graphs is that their intrinsic dimension is $D=4$. Several {\em equilibrium} characteristics of statistical physics models on planar regular random graphs are more tractable than on the regular $2D$ lattices. Non-equilibrium processes on planar regular random graphs remain unexplored. 

\bigskip\noindent
{\bf Acknowledgments.} I benefited from the correspondence with Zoltan R\'{a}cz and Sid Redner. I am grateful to the International Centre for Mathematical Sciences, Edinburgh, for the hospitality and support during the completion of this work.

\appendix
\section{Derivation of \eqref{m-av-var}--\eqref{Pmt-Phi} and \eqref{n-av-var}}
\label{ap:der}

In this Appendix, we consider the DSF process on the complete graph $K_{N+1}$ and show that the number of empty vertices $m$ is an asymptotically Gaussian random variable in the $N\to \infty$ limit. The following derivation is not fully rigorous since we assume the linear in $N$ scaling of the average and variance, $\langle m\rangle = N n(t)$ and $\langle\!\langle m^2\rangle\!\rangle = N v(t)$. We rely on the Van Kampen expansion \cite{VanKampen}, known to provide an asymptotically exact description of many reaction processes. This approach leads to a quick and efficient computation of the rescaled average $n(t)$ and variance $v(t)$ together with the rescaled distribution, Eq.~\eqref{Pmt-Phi}, of the number of empty vertices. 

The probability distribution $P_m(t)$ satisfies \eqref{Pm} with $r_m=m^2/N$, i.e.,
\begin{equation}
\label{Pmt}
N\,\frac{d P_m}{dt}=(m+1)^2 P_{m+1} - m^2 P_m
\end{equation}
We are interested in the regime where $m\sim N$, so we treat $m$ as a continuous variable and reduce a set of ordinary differential equations \eqref{Pm} into a single partial differential equation
\begin{equation}
\label{PDE}
\partial_t P=N^{-1}\left(\partial_m+\frac{1}{2}\partial_m^2\right)
\left(m^2P\right)
\end{equation}
We now transform the variables, $(m,t)\to (\xi,t)$, where $\xi = (m-Nn)/\sqrt{Nv}$, cf. \eqref{scal:def}, is a properly centralized and rescaled variable which remains finite in the $N\to\infty$ limit. The derivatives in the old and new variables are related via
\begin{subequations}
\label{derivatives}
\begin{align}
\label{der:m}
&\partial_m = \frac{\partial_\xi}{\sqrt{Nv}}\,, \quad \partial_m^2 = \frac{\partial_\xi^2}{Nv}\\
\label{der:t}
& \partial_t = \frac{\partial}{\partial t} -\left[\sqrt{N}\,\frac{\dot n}{\sqrt{v}} +\frac{\xi \dot v}{2v}\right]\partial_\xi
\end{align}
\end{subequations}
Hereinafter, we denote the time derivative of a function that depends only on time $t$ by dot. (The partial derivative $\partial_t$ in \eqref{der:t} is computed at fixed $m$, while the partial derivative $\frac{\partial}{\partial t}$ is computed at fixed $\xi$.) Substituting \eqref{derivatives} into \eqref{PDE} and using the scaling form \eqref{Pmt-scal} we obtain 
\begin{equation}
\label{eq-long}
\begin{split}
&-\sqrt{N}\,\frac{\dot n}{\sqrt{v}}\, \partial_\xi \Phi - \frac{\dot v}{2v}\,\partial_\xi(\xi \Phi) \\
&= \left(\frac{\partial_\xi}{\sqrt{Nv}}+\frac{\partial_\xi^2}{2Nv}\right)\left[\big(\sqrt{N}\,n+\xi \sqrt{v}\big)^2 \Phi\right]
\end{split}
\end{equation}
Equating the leading $O(\sqrt{N})$ terms in \eqref{eq-long},  we arrive at a `mean-field' equation for the density of empty vertices:
\begin{align}
\label{n-eq}
\dot n  =  -n^2
\end{align}
Equating the subleading $O(1)$ terms in \eqref{eq-long} yields 
\begin{equation}
\label{F-eq}
n^2\partial_\xi^2 \Phi+(\dot v+4 n v)\partial_\xi(\xi \Phi)=0
\end{equation}
Integrating \eqref{F-eq} gives $\Phi=\exp\!\left[-\frac{\dot v+4 n v}{2n^2}\xi^2\right]$ up to the normalization factor. Since the variable $\xi$ has zero mean and unit variance, $\Phi(\xi)$ must be the normal distribution \eqref{Gauss}, and therefore
\begin{align}
\label{v-eq}
\dot v + 4n v =  n^2
\end{align}
Dividing \eqref{v-eq} by \eqref{n-eq} gives $\frac{dv}{dn}=\frac{4v}{n}-1$, from which $v=Cn^4+n/3$. Solving \eqref{n-eq} and using initial conditions to fix constant $C$ we finally arrive at 
\begin{equation}
\label{n-v}
n = \frac{n_0}{1+n_0 t}\,, \quad v = \frac{n}{3} +\left(v_0-\frac{n_0}{3}\right)\left(\frac{n}{n_0}\right)^4
\end{equation}

The probability that the evolution begins with $m$ empty vertices is $\binom{N}{m}e^{-m}\big(1-e^{-1}\big)^{N-m}$ for the uncorrelated initial condition. For this binomial distribution $n_0=e^{-1}$ and $v_0=e^{-1}(1-e^{-1})$, so Eqs.~\eqref{n-v} reduce to the announced solution \eqref{m-av-var}. 

If all particles are initially at a single vertex of the complete graph, we have $n_0=1$ and $v_0=0$, so Eqs.~\eqref{n-v} reduce to the announced solution \eqref{n-av-var}. The distribution $P_m(t)$ approaches a Gaussian form, Eq.~\eqref{Gauss}, in the scaling region \eqref{scal:def} with $n$ and $v$ given by \eqref{n-av-var}. The exact solution of \eqref{Pm} subject to the initial condition $P_m(0)=\delta_{m,N}$ is also simple: The Laplace transform of $P_m(t)$ defined by \eqref{Qms:def} reads
\begin{equation}
\label{Qms}
Q_m(s)= \frac{N}{m^2}\prod_{\ell=m}^{N}  \left[1+\frac{sN}{\ell^2}\right]^{-1}
\end{equation}

The ratio of the variance of a random quanity to its average is known as a Fano factor \cite{Fano}. Equations \eqref{n-v} show that the Fano factor associated with the number of empty verticies quickly approaches a universal (independent on the initial condition) value: 
\begin{equation}
\label{Fano}
\frac{\langle\!\langle m^2\rangle\!\rangle}{\langle m\rangle} = \frac{1}{3}
\end{equation}
for $t\gg 1$. More precisely, \eqref{Fano} is valid sufficiently far from the halting time, i.e., when $1\ll t \ll N$.

\section{Higher cumulants}
\label{ap:high}
 
 Here, we present an alternative derivation of Eqs.~\eqref{n-eq} and \eqref{v-eq}. We then show how to probe higher cumulants and fully determine the third and fourth cumulants. 

Using Eqs.~\eqref{Pmt}, we deduce the governing equations for the moments $\langle m^k\rangle=\sum m^k P_m(t)$
\begin{subequations}
\label{m1234}
\begin{align}
\label{m1}
N\,\frac{d \langle m\rangle}{dt} &= - \langle m^2\rangle \\
\label{m2}
N\,\frac{d \langle m^2\rangle}{dt} &= - 2\langle m^3\rangle +  \langle m^2\rangle \\
\label{m3}
N\,\frac{d \langle m^3\rangle}{dt} &= - 3 \langle m^4\rangle +3 \langle m^3\rangle - \langle m^2\rangle \\
\label{m4}
N\,\frac{d \langle m^4\rangle}{dt} &= - 4 \langle m^5\rangle + 6 \langle m^4\rangle - 4 \langle m^3\rangle + \langle m^2\rangle 
\end{align}
\end{subequations}
etc. It is more convenient to deal with cumulants instead of the moments. We shortly write $\kappa_p =\langle\!\langle m^p\rangle\!\rangle$ for the $p^\text{th}$ cumulant. The first cumulant (equal to the first moment) appears in many following formulas, so we use the shorthand notation $\mu=\kappa_1=\langle m\rangle$. Standard formulae relate higher moments and cumulants, e.g., 
\begin{subequations}
\begin{align}
\label{cum2}
\langle m^2\rangle & = \mu^2 + \kappa_2 \\
\label{cum3}
\langle m^3\rangle & = \mu^3+3 \mu \kappa_2 + \kappa_3\\
\label{cum4}
\langle m^4\rangle &= \mu^4+6\mu^2 \kappa_2+3 \kappa_2^2 + 4\mu\kappa_3 + \kappa_4 \\
\label{cum5}
\langle m^5\rangle &= \mu^5+10\mu^3 \kappa_2+15\mu\kappa_2^2 + 10\mu^2 \kappa_3 \nonumber \\
& + 10\kappa_2 \kappa_3 + 5\mu \kappa_4  +  \kappa_5
\end{align}
\end{subequations}

Equations \eqref{m1234} for the moments are not recurrent, and hence unsolvable. However, we are interested in the asymptotic behavior when $N\gg 1$, and in this limit, the cumulants are more convenient than moments as they scale linearly in $N$ in the leading order. This assertion holds for the first cumulant: $\mu = N n$. In Appendix \ref{ap:der}, we argued in favor of the linear in $N$ scaling of the second cumulant, the variance. The same behavior continues to hold for higher cumulants
\begin{equation}
\label{rescaled}
\kappa_p = Nv_p
\end{equation}
In contrast to the moments, the rescaled cumulants satisfy recurrent, and hence solvable, equations. 

Equation \eqref{m1} becomes
\begin{equation}
\dot n = - n^2 - N^{-1}v
\end{equation}
and hence leads to \eqref{n-eq} in the $N\to\infty$ limit. Subtracting \eqref{m1} multiplied by $2\mu$ from \eqref{m2} gives the evolution equation for the variance 
\begin{equation}
\label{K2}
N\,\frac{d \kappa_2}{dt} =\mu^2 -4\mu\kappa_2 + \kappa_2 - 2 \kappa_3
\end{equation}
In deriving \eqref{K2}, we also used \eqref{cum2}--\eqref{cum3}. Recalling the linear growth, $\mu = N n$ and \eqref{rescaled}, we rewrite \eqref{K2} as 
\begin{equation}
\dot v = n^2 - 4nv +N^{-1}(v-2v_3)
\end{equation}
which reduces to \eqref{v-eq} in the $N\to\infty$ limit. 

The virtue of the above alternative derivation of \eqref{n-eq} and \eqref{v-eq} is that the same approach applies to higher cumulants. The calculations are straightforward, albeit cumbersome as the order of the cumulant increases. Below we limit ourselves to the third and fourth cumulants.  

Using Eq.~\eqref{m3} for the third moment, we derive an equation for the third cumulant 
\begin{eqnarray}
\label{K3}
N\,\frac{d \kappa_3}{dt} &=&6(\mu\kappa_2 - \mu\kappa_3 -\kappa_2^2)-\mu^2  \nonumber \\
&-& \kappa_2 + 3\kappa_3 - 3\kappa_4
\end{eqnarray}
which we combine with $\mu = N n$ and \eqref{rescaled} to find 
\begin{align}
\label{w-eq}
\dot v_3 + 6n v_3 =  -n^2+6v(n-v)
\end{align}
in the leading order. 

Dividing \eqref{w-eq} by \eqref{n-eq} gives
\begin{equation}
\frac{dv_3}{dn}=\frac{6v_3}{n}+1-\frac{6v}{n}\left(1-\frac{v}{n}\right)
\end{equation}
which has a simple polynomial solution
\begin{equation}
\label{v3}
v_3 = \frac{n}{15} + a n^4 + b n^6 + c n^7
\end{equation}
with amplitudes $a,b,c$ determined by the initial values $n(0), v(0), v_3(0)$. In the extreme case when all particles are initially at a single vertex of the complete graph, we have $n(0)=1$ and $v(0) = v_3(0) = 0$, and \eqref{v3} becomes
\begin{equation}
\label{w:extreme}
v_3 = \frac{n}{15} - \frac{n^4}{3} -\frac{2 n^6}{5}+\frac{2 n^7}{3} 
\end{equation}

Similarly, using Eq.~\eqref{m4} for the fourth moment, we derive an equation for the fourth cumulant 
\begin{eqnarray}
\label{K4}
N\,\frac{d \kappa_4}{dt} &=&\mu^2 -4\mu(2\kappa_2 - 3\kappa_3 + 2 \kappa_4)  +  12 \kappa_2(\kappa_2 -  2\kappa_3) \nonumber \\
&+& \kappa_2 - 4 \kappa_3 + 6\kappa_4 - 4 \kappa_5
\end{eqnarray}
from which we get
\begin{align}
\label{u-eq}
\dot v_4 + 8n v_4 =  n^2 - 4n(2v-3v_3)+12v(v-2v_3)
\end{align}
in the leading order. Equations \eqref{n-eq}, \eqref{v-eq}, \eqref{w-eq}, and \eqref{u-eq} for the rescaled cumulants are indeed recurrent.

Dividing \eqref{u-eq} by \eqref{n-eq} gives
\begin{equation}
\label{v4:eq}
\frac{dv_4}{dn}=\frac{8v_4}{n}-1 +4\,\frac{2v-3v_3}{n} + 12\,\frac{v(2 v_3-v)}{n^2}
\end{equation}
Recall that the rescaled cumulants $v$ and $v_3$ are the polynomials in $n$, see \eqref{n-v} and  \eqref{v3}. The general solution of Eq.~\eqref{v4:eq} is again polynomial in density: 
\begin{equation}
\label{v4}
v_4 = - \frac{n}{105}+a_4 n^4 + \sum_{j=6}^{10} a_j n^j
\end{equation}
The amplitudes in the higher-order terms depend on the initial condition. For instance
\begin{equation*}
v_4 = - \frac{n}{105} - \frac{n^4}{5} -\frac{4 n^6}{5}+\frac{4 n^7}{3}-\frac{6 n^8}{7}+\frac{16 n^9}{5} - \frac{8 n^{10}}{3}
\end{equation*}
if particles are initially at a single vertex.

The leading behavior of the cumulants is universal, i.e., independent of the initial condition when $t\gg 1$. Indeed, the linear terms in the density dominate when $n\ll 1$, and these terms are universal, cf. Eqs.~\eqref{n-v}, \eqref{v3}, \eqref{v4}. Therefore, the Fano factor \eqref{Fano} and the following two Fano factors quickly approach universal (independent of the initial condition) values
\begin{equation}
\label{Fano-34}
\frac{\kappa_2}{\kappa_1} = \frac{1}{3}\,, \quad \frac{\kappa_3}{\kappa_1} = \frac{1}{15}\,, \quad \frac{\kappa_4}{\kappa_1} = -\frac{1}{105}
\end{equation}

These Fano factors previously appeared in the context of the statistics of current in quantum conductors with many channels in the metallic regime \cite{Levitov}, quasi-classical conductors \cite{Nagaev}, and the symmetric exclusion process in one and higher dimensions \cite{Derrida04, Bodineau04, Bodineau13}. If all Fano factors associated with the random variable $m$ are the same as the Fano factors in those problems, the following Fano factors should be 
\begin{equation*}
\frac{\kappa_5}{\kappa_1} = -\frac{1}{105}\,, \quad \frac{\kappa_6}{\kappa_1} = \frac{1}{231}\,, \quad \frac{\kappa_7}{\kappa_1} = \frac{27}{5005}\,, \quad 
\frac{\kappa_8}{\kappa_1} = -\frac{3}{715}
\end{equation*}
etc. Thus, conjecturally, all Fano factors are encapsulated in the cumulant generating function 
\begin{equation}
\label{CGF}
\left[\log\!\big(\sqrt{e^z-1}+e^{z/2}\big)\right]^2=\sum_{p\geq 1}\frac{\kappa_p}{\kappa_1}\,\frac{z^p}{p!}
\end{equation}
The above pedestrian calculations of the Fano factors $\kappa_p/\kappa_1$ become unwieldy when $p$ increases. A more sophisticated approach is necessary for deriving the cumulant generating function \eqref{CGF}.

\section{DSF process on tori with localized intial distribution}
\label{ap:loc}

Consider first the DSF process on the infinite hypercubic lattices $\mathbb{Z}^d$ starting with $N\gg 1$ mobile $A$ particles in the origin. In the continuum framework applicable when $t\gg 1$, the density $a({\bf r}, t)$ satisfies the diffusion equation 
\begin{equation}
\label{DE}
\frac{\partial a}{\partial t} = D\Delta a
\end{equation}
with diffusion coefficient $D=(2d)^{-1}$ because we set the overall hopping rate of mobile particles to unity. The initial condition is $a({\bf r}, 0) = N \delta({\bf r})$.

The mobile particles occupy a ball of radius $R(t)$, so the boundary condition is
\begin{equation}
\label{a:BC}
a(r=R, t) = 0
\end{equation}
The radius grows due to the flux of the mobile particles. The corresponding Stefan-type boundary condition reads
\begin{equation}
\label{Stefan:BC}
\frac{dR}{dt} = - D\,\frac{\partial a}{\partial r}\Big|_{r=R}
\end{equation}
Although the underlying hypercubic lattice is not rotationally invariant, the density quickly becomes isotropic, $a({\bf r}, t)=a(r, t)$, as we have already tacitly assumed in writing the boundary conditions \eqref{a:BC}--\eqref{Stefan:BC}. Furthermore, the density quickly approaches the scaling form
\begin{equation}
\label{scaling}
a(r, t) = \frac{N}{R^d}\, \Phi(\xi), \qquad \xi = \frac{r}{R}
\end{equation}
Below we show that the radius increases slower than diffusively. Therefore, the left-hand side in the diffusion equation is asymtotically negligible, so \eqref{DE} simplifies to the Laplace equation $\Delta a=0$ which we combine with the scaling form \eqref{scaling} to find
\begin{equation}
\label{Lap}
\frac{d^2 \Phi}{d \xi^2} + \frac{d-1}{\xi}\,\frac{d \Phi}{d \xi} = 0
\end{equation}
The boundary condition \eqref{a:BC} becomes 
\begin{equation}
\label{Phi:BC}
\Phi(1) = 0
\end{equation}
while \eqref{Stefan:BC} reduces to
\begin{equation}
\label{Stefan}
2d R^{d+1}\,\frac{dR}{dt} = - N\,\frac{d \Phi}{d \xi}\Big|_{\xi=1}
\end{equation}
The conservation law 
\begin{equation}
\label{cons}
\int_0^R dr\,\Omega_d r^{d-1} a(r,t) + d^{-1}\Omega_d\,R^d = N
\end{equation}
asserts that the total number of existing mobile particles (the first term on the left-hand side) and annihilated mobile particles is equal to $N$. (We use standard notation $\Omega_d=2\pi^{d/2}/\Gamma(d/2)$ for the `surface area' of the unit sphere in $d$ dimensions;  $d^{-1}\Omega_d$ is the volume of the ball of unit radius in $d$ dimensions.) Substituting \eqref{scaling} into \eqref{cons} we obtain
\begin{equation}
\label{cons-xi}
N \Omega_d \int_0^1 d\xi\, \xi^{d-1} \Phi(\xi) + d^{-1}\Omega_d\,R^d = N
\end{equation}
Recalling that $N\gg 1$ and accounting that the ansatz  \eqref{scaling} is valid as long as the fraction of annihilated particles is small, we simplify \eqref{cons-xi} to  
\begin{equation}
\label{cons-int}
\Omega_d \int_0^1 d\xi\, \xi^{d-1} \Phi(\xi)  = 1
\end{equation}

In two dimensions, the solution of \eqref{Lap} subject to the boundary condition \eqref{Phi:BC} reads $\Phi=C\log \xi$. Plugging this solution into \eqref{cons-int}, we fix $C=-2/\pi$. Therefore
\begin{equation}
\label{Phi-2}
\Phi = -\frac{2}{\pi}\,\log \xi
\end{equation}
which we combine with \eqref{Stefan} to find $R = \left[\frac{2}{\pi}\,Nt\right]^\frac{1}{4}$. In higher dimensions, $d>2$, the solution of the Laplace equation \eqref{Lap} is a linear combination of two spherically symmetric solutions, the uniform and the Coulomb $\xi^{2-d}$ solutions, so the boundary condition \eqref{Phi:BC} fixes $\Phi=C[\xi^{2-d}-1]$ up to a constant. Plugging this solution into \eqref{cons-int} we fix $C=2d/[(d-2)\Omega_d]$. Therefore
\begin{subequations}
\label{Phi-R-d}
\begin{equation}
\label{Phi-d}
\Phi = \frac{2d}{(d-2)\Omega_d}\,[\xi^{2-d}-1]
\end{equation}
which we combine with \eqref{Stefan} to find the radius 
\begin{equation}
\label{R-d}
R = \left[\frac{d+2}{\Omega_d}\,Nt\right]^\frac{1}{d+2}
\end{equation}
\end{subequations}
in $d$ dimensions. In one dimension, a similar analysis gives $\Phi=1-|\xi|$ inside the droplet $-1\leq \xi \leq 1$, where $\xi=x/R$ in one dimension, and  
$R = \left[\frac{3}{2}\,Nt\right]^\frac{1}{3}$. Note that \eqref{R-d} was derived in $d>2$ dimensions, but it applies in all dimensions.

The $t^\frac{1}{d+2}$ growth of the radius of the droplet is slower than the diffusive $t^\frac{1}{2}$ growth. This feature provides an a posteriori justification of neglecting the term with time derivative in Eq.~\eqref{DE}. We also postulated that the growing droplet is asymptotically spherical. The droplet's shape in a similar internal diffusion-limited aggregation (DLA) problem has been investigated \cite{Bramson-IDLA, Levine12, Asselah13}. In the IDLA, the infinite hypercubic lattice $\mathbb{Z}^d$ is initially empty (equivalently, occupied by immobile $B$ particles, one particle per site), and mobile $A$ particles are injected into the origin at such a slow rate that each $A$ particle undergoes reaction before the next injection event. The droplet is asymptotically a ball, and deviations from the spherical shape are small, see \cite{Bramson-IDLA, Levine12, Asselah13}. The droplet in our problem should also be very close to spherical. 

The above analysis gives the leading asymptotic behavior of the DSF process on the {\em infinite} lattices $\mathbb{Z}^d$. The density is compact and hence describes the DSF on tori $(R_L)^d$ with $N=L^d$ vertices when $R<\frac{1}{2}L$. Comparing \eqref{R-d} with $L=N^\frac{1}{d}$, we find that \eqref{scaling} is applicable up to time $N^\frac{2}{d}$, providing the lower bound for the halting time. In one dimension, the density at the origin is $a=N/R$, hence of the order of unity when $R\sim N$. Thus, the halting time is $T_N\sim N^2$ in one dimension. In $d>2$ dimensions, the density near the origin [cf. \eqref{scaling} and \eqref{Phi-d}] is $a\sim \frac{N}{R^d} \left(\frac{1}{R}\right)^{2-d}\sim \frac{N}{R^2}$. Hence, when the droplet becomes comparable to the system size, $R\sim N^\frac{1}{d}$, the density at the origin is still very large, $a\sim  N^\frac{d-2}{d}$. We thus expect the halting time to grow faster than $N^\frac{2}{d}$ when $d>2$.

\bibliography{references-packing}

\end{document}